\documentclass[12pt,letterpaper,twoside]{article}

\usepackage{amsmath,amssymb}
\usepackage{amsthm}
\usepackage{mathtools}
\usepackage{stmaryrd}
\usepackage{fancyhdr}
\usepackage{times}
\usepackage{mathrsfs} 
\usepackage{titlesec}
\usepackage[normalem]{ulem}
\usepackage[pdftex,dvips]{geometry}

\mathtoolsset{showonlyrefs,centercolon,mathic}

\newtheoremstyle{fact}
     {\topsep}
     {\topsep}
     {\slshape}
     {}
     {\bfseries}
     {}
     { }
     {\thmname{#1}\thmnumber{ #2.}\thmnote{ \rm (#3)}}

\newtheoremstyle{mylabel}
     {\topsep}
     {\topsep}
     {\itshape}
     {}
     {\bfseries}
     {}
     { }
     {\thmname{#1}\thmnote{ #3}.}

\newtheorem{theorem}{Theorem}[section]
\newtheorem{Ltheorem}{Theorem}

\newtheorem*{theorem*}{Theorem} 
\newtheorem{lemma}[theorem]{Lemma}
\newtheorem{proposition}[theorem]{Proposition}
\newtheorem{corollary}[theorem]{Corollary}
\newtheorem{problem}{Problem}

\newtheorem*{conjecture*}{Conjecture}

\theoremstyle{mylabel}
\newtheorem*{Ltheorem*}{Theorem}

\theoremstyle{definition}
\newtheorem{definition}[theorem]{Definition}

\newtheorem*{remark*}{Remark}

\newtheorem*{question*}{Question}
\newtheorem*{examples*}{Examples}  
\newtheorem{example}[theorem]{Example}
\newtheorem{examples}[theorem]{Examples}
\newtheorem*{example*}{Example}

\newtheorem*{convention*}{Convention}

\theoremstyle{fact}

\newtheorem{ftheorem}[theorem]{Theorem}
\newtheorem{flemma}[theorem]{Lemma}

\newtheorem{fcorollary}[theorem]{Corollary}

\newenvironment{myromanlist}[1][enumi]{\begin{list}{{\rm (\roman{#1})}}
{\usecounter{#1}\setlength{\labelwidth}{25pt}\setlength{\topsep}{-6pt}
\setlength{\itemsep}{-4pt} \setlength{\leftmargin}{25pt}}}{\end{list}}

\newenvironment{myalphlist}[1][enumi]{\begin{list}{{\rm (\alph{#1})}}
{\usecounter{#1}\setlength{\labelwidth}{25pt}\setlength{\topsep}{-6pt}
\setlength{\itemsep}{-4pt} \setlength{\leftmargin}{25pt}}}{\end{list}}

\newenvironment{mynumlist}[1][enumii]{\begin{list}{{\rm (\arabic{#1})}}
{\usecounter{#1}\setlength{\labelwidth}{25pt}\setlength{\topsep}{-6pt}
\setlength{\itemsep}{-4pt} \setlength{\leftmargin}{25pt}}}{\end{list}}

\def\proofont{\fontseries{bx}\fontshape{sc}\selectfont}
\def\proofname{Proof.}

\newcommand{\Note}[1]{}

\makeatletter
\renewenvironment{proof}[1][\proofname]{\par
  \normalfont
  \topsep6\p@\@plus6\p@ \trivlist
  \item[\hskip\labelsep\noindent\proofont #1]\ignorespaces
}{%
  \qed\endtrivlist
}
\makeatother

\titlelabel{\thetitle.\ }
\titleformat*{\section}{\normalsize\bfseries\centering}
\titleformat*{\subsection}{\normalsize\bfseries\itshape}

\author{D. Dikranjan\thanks{The first author acknowledges the financial 
aid received from ``Progetti di Eccellenza 2011/12'' of Fondazione 
CARIPARO.} { }and G\'abor Luk\'acs\thanks{The second author gratefully 
acknowledges the generous financial support received from NSERC, which 
enabled him to do this research.}}

\title{Compact-like abelian groups without\\ non-trivial 
quasi-convex  null sequences}

\begin{document}

\makeatletter
\def\@fnsymbol#1{\ifcase#1\or * \or 1 \or 2  \else\@ctrerr\fi\relax}
\let\mytitle\@title

\chead{\small\itshape D. Dikranjan and G. Luk\'acs / Compact-like abelian 
groups without quasi-convex null sequences}
\fancyhead[RO,LE]{\small \thepage}
\makeatother

\maketitle

\def\thanks#1{} 

\thispagestyle{fancy}

\begin{abstract} 

In this paper, we study precompact abelian groups $G$ that contain no 
sequence $\{x_n\}_{n=0}^\infty$ such that $\{0\} \cup \{\pm x_n \mid n \in 
\mathbb{N}\}$ is infinite and quasi-convex in $G$, and $x_n 
\longrightarrow 0$. We characterize groups with this property in the 
following classes of  groups:

\begin{myalphlist}

\item
bounded precompact abelian groups;

\item
minimal abelian groups;

\item
totally minimal abelian groups; 

\item
$\omega$-bounded abelian groups.

\end{myalphlist}

\vspace{3pt}

\noindent We also provide examples of minimal abelian groups with this 
property, and show that there exists a minimal pseudocompact abelian group 
with the same property; furthermore, under Martin's Axiom, the group may 
be chosen to be countably compact minimal abelian.

\end{abstract}

\section{Introduction}

\thispagestyle{empty}

This note is a sequel to \cite{DikGL2}, and its aim is to provide an  answer to Problem I posed in that paper.

One of the main sources of inspiration for the theory of topological 
groups is the theory of topological vector spaces, where the notion of 
convexity plays a prominent role. In this context, the reals $\mathbb{R}$ 
are replaced with the circle group $\mathbb{T} := \mathbb{R}/\mathbb{Z}$, 
and linear functionals are replaced by {\em characters}, that is, 
continuous homomorphisms to $\mathbb{T}$. By making substantial use of 
characters, Vilenkin introduced the notion of quasi-convexity for abelian 
topological groups as a counterpart of convexity in topological vector 
spaces (cf.  \cite{Vilenkin}).

Let $\pi\colon \mathbb{R} \rightarrow \mathbb{T}$  denote  the  canonical projection. Since the restriction 
$\pi_{|(-\frac 1 2,\frac 1 2]}\colon (- \frac 1 2,\frac 1 2] \rightarrow \mathbb{T}$ 
is a~bijection,  we often identify in what follows, {\em par abus de language},  a~number 
$a\in (-\frac 1 2,\frac 1 2]$ with its image (coset) $\pi(a) =  a+\mathbb{Z} \in \mathbb{T}$.  
We put $\mathbb{T}_m  :=\pi([-\frac{1}{4m},\frac{1}{4m}])$ for all $m\in \mathbb{N}\backslash\{0\}$.  According to standard 
notation in this area, we use $\mathbb{T}_+$ to denote $\mathbb{T}_1$. For an abelian topological group $G$, we denote by 
$\widehat{G}$  the {\em Pontryagin dual} of $G$, that is, the group of all  characters of $G$ endowed with the compact-open topology.

\begin{samepage}
\begin{definition}\label{def:intro:qc} For $E  \subseteq G$ and 
$A \subseteq \widehat{G}$, the {\em polar} of $E$ and the {\em prepolar}  $A$ are defined as
\begin{align*} 
E^\triangleright=\{\chi\in \widehat{G} \mid \chi(E) \subseteq 
\mathbb{T}_+\} \quad \text{and} \quad
A^\triangleleft=\{ x \in G \mid  \forall \chi \in A, 
\chi(x) \in \mathbb{T}_+ \}.
\end{align*}
The set $E$ is said to be {\em quasi-convex} if  $E  = E^{\triangleright\triangleleft} $. 
\end{definition}
\end{samepage}

Obviously,  $E\subseteq  E^{\triangleright\triangleleft}$ holds for  every  $E \subseteq G$. 
Thus, $E$ is quasi-convex if and only if for every  $x\in G\backslash E$ there exists 
$\chi \in E^\triangleright$ such that  $\chi(x) \notin \mathbb{T}_+$. The set  
$Q_{G}(E)  := E^{\triangleright\triangleleft}$ is the smallest quasi-convex set of $G$ that contains $E$, and it is  called the {\em quasi-convex hull} of $E$. 

\begin{definition}
A sequence  $\{x_n\}_{n=0}^\infty \subseteq G$  is said to be {\em  quasi-convex} if  
$S  =  \{0\} \cup\{\pm x_n \mid n \in \mathbb{N}\}$ 
is quasi-convex in $G$. We say that $\{x_n\}_{n=0}^\infty$ is 
{\em non-trivial} if the set $S$ is infinite, and it is a {\em null 
sequence} if $x_n \longrightarrow 0$.
\end{definition}

\begin{examples}  \label{exp:intro:seqs}
It turns out that most ``common" groups contain a non-trivial 
quasi-convex null sequence:

\begin{myalphlist}

\item
For every prime $p$, the group $\mathbb{J}_p$ of $p$-adic integers 
admits a  non-trivial quasi-convex null sequence contained in
the subgroup $\mathbb{Z}$ (cf.~\cite[1.4]{DikLeo}, 
\cite[Theorem~D]{DikGL1}, and \cite[Theorem~B]{DikGL2}).

\item
For every prime $p$, $\mathbb{T}$ admits a  non-trivial quasi-convex  null sequence contained in the  $p\mbox{-}$com\-ponent of $\mathbb{T}$
(cf.~\cite[1.1]{DikLeo}, \cite[Theorem~B]{DikGL1}, and \cite[Theorem~C]{DikGL2}).

\item
If $\{m_k\}_{k=1}^\infty$ is a sequence of integers such that  $m_k \geq  4$ for every  $k \in \mathbb{N}$, then $\prod\limits_{k=0}^\infty \mathbb{Z}_{m_k}$ admits
a~non-trivial  quasi-convex null sequence contained in the subgroup $\bigoplus\limits_{k=0}^\infty \mathbb{Z}_{m_k}$  (cf.~\cite[5.5]{DikGL2}).
\end{myalphlist}
\end{examples}

In our paper \cite{DikGL2}, we characterized the locally compact abelian  groups that admit no non\mbox{-}trivial quasi-convex null sequences as  follows.

\begin{ftheorem}[{\cite[Theorem~A]{DikGL2}}]\label{thm:intro:LCAqcs}
For every locally compact abelian group $G$, the following statements are  equivalent:

\begin{myromanlist}

\item
$G$ admits no non-trivial quasi-convex null sequences;

\item
one of the subgroups $G[2] = \{g \in G \mid 2g = 0\}$ or 
$G[3] = \{g \in G \mid 3g = 0\}$ is open in $G$;

\item
$G$ contains an open compact subgroup of the form $\mathbb{Z}_2^\kappa$ or 
$\mathbb{Z}_3^\kappa$ \ for some cardinal $\kappa$.

{\parindent -25pt 
Furthermore, if $G$ is compact, then these conditions  are also 
equivalent to:
}

\item
$G  \cong \mathbb{Z}_2^\kappa  \times  F$ or 
$G \cong  \mathbb{Z}_3^\kappa  \times  F$,
where $\kappa$ is some cardinal and $F$ is a finite abelian group;

\item
one of the subgroups $2G$ and $3G$ is finite.

\end{myromanlist}
\end{ftheorem}

We also asked whether it was possible to replace the class of locally 
compact abelian groups with a different class that contains all compact 
abelian groups (cf.~\cite[Theorem~A]{DikGL2}). In this note,~we present 
several classes, and characterizations of groups in these classes that 
admit no non-trivial quasi-convex null sequences. To that end, we recall 
some compactness-like properties.

\begin{definition} Let $G$ be a (Hausdorff) topological group.

\begin{myalphlist}

\item
$G$ is {\itshape precompact} if it can be covered by finitely many 
translates of any neighborhood of the identity, or equivalently, if it is 
a dense subgroup of a compact group $\widetilde G$, its {\itshape completion};

\item
$G$ is {\itshape minimal} if there is no coarser (Hausdorff) group 
topology (cf.~\cite{Steph} and~\cite{Doi});

\item
$G$ is {\itshape totally minimal} if every (Hausdorff) quotient group $G$ 
is minimal (cf.~ \cite{DikPro});

\item $G$ is {\itshape pseudocompact} if every real-valued continuous 
function on $G$ is bounded (cf.~ \cite[3.10]{Engel6});

\item
$G$ is {\itshape countably compact} if every countable open cover of $G$ 
admits a finite subcover (cf.~ \cite[3.10]{Engel6});

\item
$G$ is {\itshape $\omega$-bounded} if every countable subset of $G$ is  contained in a compact subgroup of $G$.

\end{myalphlist}
\end{definition}

By the celebrated Prodanov-Stoyanov Theorem, every minimal abelian group  is precompact (cf.~\cite{ProdStoj} and~\cite{ProdStoj2}), and thus the 
relationships among the aforementioned properties can be described as  follows.
\begin{align}
\mbox{$\omega$-bounded} \Longrightarrow
\mbox{countably compact} \Longrightarrow
\mbox{pseudocompact}  \Longrightarrow
\mbox{precompact} \\
\mbox{totally minimal} \Longrightarrow
\mbox{minimal}  \stackrel{\text{abelian}}\Longrightarrow
\mbox{precompact} 
\end{align}
For greater clarity, all of these implications, except for the very last 
one, hold without the assumption that the group in question is abelian.

The classes of groups studied in this paper overlap with those 
investigated in \cite{DikGL2} only to the smallest possible extent, 
because every group that is both precompact and locally compact is 
actually compact. Thus, the present note is complementary to our work in 
\cite{DikGL2}. Our first two results demonstrate the level of complexity 
of the problem of finding non-trivial quasi-convex null sequences when one 
leaves the class of locally compact abelian groups. Indeed, in the absence 
of local compactness, Theorem~\ref{thm:intro:LCAqcs} may fail even in the 
presence of strong compactness-like properties.

Recall that the {\itshape exponent} of an abelian group $G$ is 
$\exp G := \inf\{ n >0 \mid nG=0\}$. The group $G$ is {\itshape bounded} 
if it has a finite exponent.

\begin{Ltheorem} \label{thm:main:ABSbad}
Let $p = 2$ or $p = 3$, and let $\kappa$ be an infinite cardinal.

\begin{myalphlist}

\item
There  exists a minimal abelian  group $G$ of exponent $p^2$ such that $|pG| = \kappa$ and $G$ admits no non-trivial quasi-convex null sequences.

\item
If $\kappa^\omega=\kappa$, then there exists a  minimal  \uline{pseudocompact} abelian group of exponent $p^2$ such  that 
$|pG|=\kappa$ and $G$ admits no non-trivial  quasi-convex null sequences.

\end{myalphlist}
\end{Ltheorem}

In part (b) of Theorem~\ref{thm:main:ABSbad}, $pG$ is a pseudocompact 
group (being a continuous image of $G$),~and thus $|pG|=\kappa$ must 
satisfy certain constraints that the size of every infinite pseudocompact 
homogeneous space does: $\kappa \geq \mathfrak{c}$, and $\kappa$ cannot be 
a strong limit cardinal of countable cofinality (cf.~\cite[1.2, 
1.3(a)]{vanDouwPS}). We note that both of these conditions follow from the 
hypothesis $\kappa^\omega=\kappa$. Under the Generalized Continuum 
Hypothesis (GCH), $\kappa^\omega > \kappa$ for a cardinal 
$\kappa\geq\mathfrak{c}$ if and only if $\kappa$ is a strong limit 
cardinal of countable cofinality. Therefore, under GCH, the hypothesis 
$\kappa^\omega=\kappa$ is not only sufficient but also {\itshape 
necessary} for the existence of a group $G$ as in 
Theorem~\ref{thm:main:ABSbad}(b).

\begin{Ltheorem} \label{thm:main:MACHbad}
Under Martin's Axiom (or the Continuum Hypothesis),
there exists a countably compact minimal abelian group of 
exponent $4$ such that $2G$ is infinite and $G$ admits no non-trivial 
quasi-convex null sequences.
\end{Ltheorem}

The proofs of Theorems~\ref{thm:main:ABSbad} and~\ref{thm:main:MACHbad} 
are presented in \S\ref{sect:bad}, and they are based on ``lifting" known 
examples of precompact abelian groups of exponent $p$ that admit no 
non-trivial convergent sequences at all (cf.~\cite{Flor}, 
\cite[Theorem~3]{DijvMil}, \cite[8.1 \& 9]{vDouw}, and~\cite{HajJuhCC}). 
Although  Theorem~\ref{thm:main:MACHbad}  may appear to suggest a  
negative answer to \cite[Problem~II]{DikGL2}, it is possible to establish 
a criterion in the spirit of Theorem~\ref{thm:intro:LCAqcs} for 
precompact groups of a finite exponent. To that end, recall that
a subset $A$ of a topological space $X$ is {\itshape sequentially open} if 
for every convergent sequence $\{x_n\}\subseteq X$ such that $\lim x_n \in 
A$, one has $x_n \in A$ for all but finitely many $n$.

\pagebreak[3]

\begin{samepage}
\begin{Ltheorem} \label{thm:main:BP}
For every bounded precompact abelian group $G$, the following statements 
are equivalent:

\begin{myromanlist}

\item
$G$ admits no non-trivial quasi-convex null sequences;

\item
one of the subgroups  $G[2]  = \{g \in G \mid 
2g =  0\}$ or $G[3]  = \{g  \in  G \mid  3g  =  0\}$ is sequentially open 
in~$G$.

\end{myromanlist}

\end{Ltheorem}
\end{samepage}


The proof of Theorem~\ref{thm:main:BP} is presented in \S\ref{sect:BP}.

Certainly, the entire problem of searching for non-trivial quasi-convex 
null sequences is meaningless in a group that admits no non-trivial 
convergent sequences at all. While $\omega$-bounded groups are known to 
have non-trivial convergent sequences, our interest in minimal abelian 
groups is also motivated by a recent result of Shakhmatov, which 
guarantees the existence of non-trivial convergent sequences in minimal 
abelian groups (cf.~\cite[1.3]{Shakh-MinAbSeq}).

\begin{samepage}
\begin{Ltheorem} \label{thm:main:minimal}

For every minimal abelian group $G$, the following statements are  equivalent:

\begin{myromanlist}
\item
$G$ admits no non-trivial quasi-convex null sequences;

\item
$G \cong P \times F$, where $P$ 
is a minimal bounded abelian $p$-group {\rm (}$p \leq 3${\rm )} admitting 
no non-trivial quasi-convex null sequences, and $F$ is a finite abelian 
group;

\item
one of the subgroups  $G[2]  = \{g \in G \mid  2g =  0\}$ or 
$G[3]  = \{g  \in  G \mid  3g  =  0\}$ is sequentially open 
in~$G$;

\item
$G$ contains a sequentially open compact subgroup of the form 
$\mathbb{Z}_2^\kappa$ or $\mathbb{Z}_3^\kappa$ \ for some cardinal $\kappa$.

\end{myromanlist}

\end{Ltheorem}
\end{samepage}

The following result shows that much of Theorem~\ref{thm:intro:LCAqcs} can 
be salvaged by imposing stronger compactness-like properties.

\begin{Ltheorem} \label{thm:main:omega}
The following statements are equivalent for every abelian group $G$ that 
is  $\omega$-bounded  \underline{or} totally minimal:

\begin{myromanlist}

\item
$G$ admits no non-trivial quasi-convex null sequences;

\item
one of the subgroups $G[2]  = \{g \in G \mid 2g =  0\}$ or
$G[3]  = \{g  \in  G \mid 3g  =  0\}$ is open in $G$;

\item
one of the subgroups $2G$ and $3G$ is finite.

{\parindent -25pt 
Furthermore, if $G$ is totally minimal, then these conditions are also 
equivalent to:
}

\item
$G  \cong \mathbb{Z}_2^\kappa  \times  F$ or  
$G \cong  \mathbb{Z}_3^\kappa  \times  F$, where $\kappa$ is some 
cardinal and $F$ is a finite abelian group.

\end{myromanlist}
\end{Ltheorem}

Theorems \ref{thm:main:ABSbad} and \ref{thm:main:MACHbad} show that in 
Theorem~\ref{thm:main:omega}, one 
cannot weaken  ``$\omega$\mbox{-}bounded" to  ``countably compact and 
minimal" (or to ``pseudocompact and minimal"), and ``totally minimal" to 
``minimal" in Theorem~\ref{thm:main:omega}. Therefore,
Theorem~\ref{thm:main:minimal} is the best one can achieve in the class 
of minimal abelian groups.  The proofs of Theorems~\ref{thm:main:minimal} 
and~\ref{thm:main:omega}  are  presented in \S\ref{sect:omega}. The 
totally minimal case of  Theorem~\ref{thm:main:omega} relies on 
intermediate steps in the proof of  Theorem~\ref{thm:main:minimal}. One of 
the main ingredients of Theorem~\ref{thm:main:minimal} is the following 
result.

\begin{Ltheorem} \label{thm:main:8T}
Let $\{q_n\}_{n=0}^\infty$ be a sequence of positive integers, and put 
$b_n = q_0\cdots q_n$ \ for every $n \in \mathbb{N}$. If $q_n \geq 8$ for 
every $n \in \mathbb{N}$, then $\{\frac 1 {b_n} \}_{n=0}^\infty$ is a 
quasi-convex sequence in $\mathbb{T}$.
\end{Ltheorem}

Theorem~\ref{thm:main:8T} implies that the subgroup of 
$\mathbb{Q}/\mathbb{Z}$ generated by elements of prime order admits a 
non-trivial quasi-convex null sequence. We note that the condition $q_n 
\geq 8$ in Theorem~\ref{thm:main:8T} is unnecessarily restrictive, and can 
be replaced with $q_n \geq 5$; however, the proof of the latter is more 
complicated and longer, and Theorem~\ref{thm:main:8T} is sufficient for 
establishing Theorem~\ref{thm:main:minimal}. The proof of 
Theorem~\ref{thm:main:8T} is presented in \S\ref{sect:T}.

\pagebreak[3]

Finally, we turn to posing some open problems. As we mentioned earlier, 
under GCH, the hypothesis $\kappa^\omega=\kappa$ is not only sufficient 
but also necessary for the existence of a group $G$ as in 
Theorem~\ref{thm:main:ABSbad}(b). Thus, in ZFC, it is not possible to 
prove Theorem~\ref{thm:main:ABSbad}(b) without the 
assumption~$\kappa^\omega=\kappa$.

\begin{problem}
Let $p=2$ or $p=3$, and let
$\kappa$ be a cardinal such that $\kappa \geq \mathfrak{c}$ and 
$\kappa^\omega > \kappa$. Is it consistent that there exists a~minimal 
\uline{pseudocompact} abelian group of exponent $p^2$ such that 
$|pG|=\kappa$ and $G$ admits no non-trivial quasi-convex null sequences?
\end{problem}

We show in \S\ref{sect:bad} that if one of the subgroups $G[2]$~or $G[3]$ 
is sequentially open in a (locally) precompact group $G$, then it contains 
no non-trivial quasi-convex null sequences 
(Proposition~\ref{prop:bad:G23}).  Theorems~\ref{thm:main:BP}, 
\ref{thm:main:minimal}, and \ref{thm:main:omega} state that the converse 
of this implication is also true in each of the classes of bounded 
precompact, minimal, totally minimal, and $\omega$-bounded abelian groups. 
Thus, it is natural to ask whether the implication remains reversible in 
general, without assuming some compact-like properties.

\begin{problem}
Let $G$ be a (locally) precompact group, and suppose that $G$ admits no 
non-trivial quasi-convex null sequences. Is  one of the subgroups $G[2]$  
or $G[3]$ sequentially open in $G$?
\end{problem}

\section{Preliminaries}

\label{sect:prel}

In this section, we provide a few well-known definitions and results that we rely on in the paper. We start off by recalling some 
terminology from duality theory. Let $H$ be a subgroup of an abelian topological group $G$. The 
{\itshape annihilator} of $H$ in $\widehat G$ is the subgroup $H^\perp  :=  \{\chi \in \widehat G \mid\chi(H)= \{0\}\}$. 
The subgroup $H$ is {\itshape dually closed} in $G$ if $H =  \bigcap\{\ker \chi \mid 
\chi  \in H^\perp\}$. Since $H^\perp  =  H^\triangleright$ for every 
subgroup, $H$ is dually closed in $G$ if and only if it is quasi-convex in $G$.
The subgroup $H$ is {\itshape dually embedded} in $G$ if every continuous character 
of $H$ has an extension to a continuous character of~$G$, that is, the restriction homomorphism $\widehat G \rightarrow \widehat H$ is surjective.

\begin{examples} \mbox{ } \label{exp:prel:dual}
\begin{myalphlist}

\item
If $H$ is an open subgroup of the abelian topological group $G$, then $H$ is  dually closed and dually embedded in $G$ (cf.~\cite[3.3]{Noble}).

\pagebreak[3]

\item
If $H$ is a closed subgroup of a locally compact abelian group $L$, then 
$H$ is dually closed and dually embedded in $L$ (cf.~\cite[Theorems~37 
and~54]{Pontr}).

\item
If $H$ is a dense subgroup of an abelian topological group $G$, then $H$ 
is dually embedded in $G$.


\item
Every subgroup of a locally precompact abelian group is dually embedded in 
it.

\end{myalphlist} 
\end{examples}

\begin{lemma} \label{lemma:prel:subgr}
Let $G$ be an abelian topological group. 

\begin{myalphlist}

\item
If $H$ is a dually embedded subgroup of $G$, then $Q_H(S)  = Q_G(S)\cap H$ 
for every subset~$S$~of~$H$.

\item
If $H$ is a dually closed and dually embedded subgroup of $G$, then $Q_H(S)  =  Q_G(S)$ for every subset $S$ of~$H$.

\end{myalphlist}
\end{lemma}

Lemma~\ref{lemma:prel:subgr} is similar to \cite[5.1]{DikGL2}, and its proof, which relies on the following general property of the quasi-convex 
hull, is provided here only for the sake of completeness.

\begin{flemma}[{\cite[I.3(e)]{GLdualtheo}, \cite[2.7]{DikLeo}}]
\label{lemma:prel:qc-hom}
If $f\colon G\to H$ is a continuous homomorphism of abelian topological groups, and $E\subseteq G$, then $f(Q_G(E)) \subseteq Q_H(f(E))$.
\end{flemma}

\begin{proof}[Proof of Lemma~\ref{lemma:prel:subgr}.]
Let $S \subseteq H$ be a subset.

(a) Let $\iota\colon H \rightarrow G$ denote the inclusion map.By Lemma~\ref{lemma:prel:qc-hom},
\begin{align}
Q_H(S) \subseteq \iota^{-1}(Q_G(S))=Q_G(S)\cap H.
\end{align}
To show the reverse inclusion, let 
$h \in 
Q_G(S) \cap  H $. If 
$\chi \in \widehat H$ is such that
$\chi(S) \subseteq \mathbb{T}_+$, then 
$\chi= \psi_{|H}$ for some 
$\psi \in \widehat G$ (as $H$ is dually embedded in $G$), and one has $\psi(S) = \chi(S) 
 \subseteq  \mathbb{T}_+$. Consequently,
$\chi(h) = \psi(h) \in \mathbb{T}_+$, 
because $h \in Q_G(S)$. This shows that
$Q_G(S)\cap H \subseteq Q_H(S)$, as  required.

(b) Since $H$ is dually closed, $Q_G(S)\subseteq Q_G(H)  =  H$. Consequently, the statement follows from part~(a).
\end{proof}

Examples~\ref{exp:prel:dual} combined with Lemma~\ref{lemma:prel:subgr}  yields the following consequences.

\begin{corollary} \label{cor:prel:subgrpB}
Let $G$ be a {\rm (}locally{\rm )} precompact abelian group, and $H$ a 
closed subgroup. Then $H$ is dually closed and dually embedded in $G$, and 
if $\{x_n\} \subseteq H$ is a quasi-convex sequence, then $\{x_n\}$ is 
quasi-convex in $G $.
\end{corollary}

\begin{proof}
By Example~\ref{exp:prel:dual}(d), $H$ is dually 
embedded in $G$. By Example~\ref{exp:prel:dual}(b), 
$\operatorname{cl}_{\widetilde{G}} H$ is dually closed in the 
completion $\widetilde G$. Thus, by Lemma~\ref{lemma:prel:subgr}(a),
\begin{align}
 Q_G(H) = Q_{\widetilde G}(H) \cap G \subseteq 
Q_{\widetilde G}(\operatorname{cl}_{\widetilde G} H) \cap G \subseteq 
(\operatorname{cl}_{\widetilde G} H)\cap G = H, 
\end{align}
because $H$ is closed in $G$. This shows that $H$ is dually closed in 
$G$. The second statement follows now by Lemma~\ref{lemma:prel:subgr}(b).
\end{proof}

\begin{corollary} \label{cor:prel:subgrpA}
Let $G$ be a {\rm (}locally{\rm )} precompact abelian group, $H$ a (not 
necessarily closed) subgroup, $\{x_n\} \subseteq G$ a quasi-convex 
sequence such that  $x_n \in H$ for infinitely many  $n \in\mathbb{N}$, 
and $\{x_{n_k}\}$ the subsequence of $\{x_n\}$ consisting of all members 
that belong to $H$. Then $\{x_{n_k}\}$ is quasi-convex in $H$.
\end{corollary}

\begin{proof} 
By Example~\ref{exp:prel:dual}(d), $H$ is dually embedded in $G$.~Put 
$S  :=  \{\pm x_n \mid n  \in \mathbb{N}\}  \cup  \{0\}$ and
set $S^\prime   := \{\pm x_{n_k} \mid k  \in  \mathbb{N}\} \cup  \{0\}$. By 
Lemma~\ref{lemma:prel:subgr}(a), 
\begin{align}
Q_H(S^\prime) = Q_G(S^\prime) \cap H \subseteq Q_G(S) \cap H = S \cap H. 
\end{align}
It follows from the definition of the subsequence $\{x_{n_k}\}$ that 
$S \cap H =  S^\prime $, as desired.
\end{proof}

We turn now to minimality and total minimality.

\begin{ftheorem} \label{thm:prel:ess} 
Let $G$ be an abelian group with completion $\widetilde G$. Then:

\begin{myalphlist}

\item
{\rm ({\cite[Theorem~2]{Steph},\cite{Prod1}, 
\cite[Propositions~1~and~2]{Banasch}, \cite[3.31]{GLCLTG}})}
$G$ is minimal if and only if $G$ is  precompact and $G \cap H \neq \{0\}$ for every non-trivial closed subgroup $H$ of $\widetilde G$;

\item
{\rm (\cite{DikPro}, \cite[3.31]{GLCLTG})}
$G$ is totally minimal if and only if $G$ is precompact and $G\cap H$ is dense in $H$ for every closed subgroup of $\widetilde G$.

\end{myalphlist}

\end{ftheorem}

Recall that the {\itshape socle} of $\operatorname{soc}(A)$ of an abelian 
group $A$ is the subgroup of torsion elements whose order is square-free
(that is, not divisible by the square of a prime number), 
or equivalently, the direct sum of the subgroups 
$A[p] :=\{ x\in A \mid px = 0\}$, where $p$ ranges over all primes. We 
put $\operatorname{tor}(A)$ for the torsion subgroup of $A$.

\begin{fcorollary} \label{cor:prel:st}
Let $G$ be an abelian group with completion $\widetilde G$.
\begin{myalphlist}

\item
If $G$ is minimal, then
$\operatorname{soc}(\widetilde G)
\subseteq G $.

\item
{\rm (\cite[4.3.4]{DikProSto})}
If $G$ is totally minimal, then
$\operatorname{tor}(\widetilde G)
\subseteq  G $.

\end{myalphlist}

\vspace{6pt}

\noindent
Furthermore, if $\widetilde G$ is a bounded compact abelian group, then the converse of {\rm (a)} is also true.
\end{fcorollary}

\begin{definition} (\cite{DikProET}, \cite[p.~141]{DikProSto})
A compact abelian group is an {\em exotic torus} if it contains no 
subgroup that is topologically isomorphic to the $\mathbb{J}_p$ 
($p$-adics) for some prime $p$.
\end{definition}

The notion of exotic torus was introduced by Dikranjan and Prodanov in 
\cite{DikProET}, who also provided, among other things, the following 
characterization for such groups.

\begin{ftheorem}[{\cite{DikProET}}, {\cite[2.6]{DikGL2}}] \label{thm:prel:ET}
A compact abelian group $K$ is an exotic torus if and only if it contains a~closed subgroup $B$ such that

\begin{myromanlist}

\item
$K/B \cong \mathbb{T}^n$ for some $n \in \mathbb{N}$, and

\item
$B = \prod\limits_{p} B_p$, where each $B_p$ is a compact bounded abelian  $p$-group.

\vspace{6pt}

\end{myromanlist}
Furthermore, if $K$ is connected, then each $B_p$ is finite.
\end{ftheorem}

Recall that a topological group is {\em pro-finite} if it is the  (projective) limit of finite groups, or equivalently, if it is compact and 
zero-dimensional. For a prime $p$, a topological group $G$ is called a {\em pro-$p$-group} if it is the (projective) limit of finite 
$p$-groups, or equivalently, if it is pro-finite and $x^{p^n} \longrightarrow e$ for every $x \in G$ (or,  in the abelian case, $p^n x \longrightarrow 0$).

\begin{ftheorem}[{\cite{ArmacostLCA}, 
\cite[Corollary 8.8(ii)]{HofMor}, \cite[4.1.3]{DikProSto}}]
\label{thm:prel:profin}
Let $G$ be an abelian pro-finite group. Then $G = \prod\limits_{p}  G_p$, 
where each $G_p$ is a pro-$p$-group.
\end{ftheorem}

Finally, we note for the sake of clarity that in our notation, 
$\mathbb{N}=\{0,1,2,\ldots\}$, that is, $0 \in \mathbb{N}$.

\section{Counterexamples}

\label{sect:bad}

\begin{Ltheorem*}[\ref{thm:main:ABSbad}]
Let $p = 2$ or $p = 3$, and let $\kappa$ be an infinite cardinal.

\begin{myalphlist}

\item
There  exists a minimal abelian  group $G$ of exponent  $p^2$ such that
$|pG| = \kappa$ and $G$ admits no non-trivial quasi-convex 
null sequences.

\item
If $\kappa^\omega=\kappa$, then there exists a  minimal
\uline{pseudocompact} abelian group of exponent $p^2$ such  that
$|pG|=\kappa$ and $G$ admits no non-trivial  quasi-convex 
null sequences.

\end{myalphlist}
\end{Ltheorem*}

\begin{Ltheorem*}[\ref{thm:main:MACHbad}]
Under Martin's Axiom (or the Continuum Hypothesis), there exists a 
countably compact minimal abelian group of exponent $4$ such that $2G$ is 
infinite and $G$ admits no non-trivial quasi-convex null sequences.
\end{Ltheorem*}

In this section, we prove Theorems~\ref{thm:main:ABSbad} 
and~\ref{thm:main:MACHbad}.

\begin{flemma}[{\cite[5.3]{DikGL2}}] \label{lemma:bad:Z23}
If $G$ is an abelian topological group of exponent $2$ or $3$, then $G$ admits no non-trivial quasi-convex null sequences.
\end{flemma}

\begin{proposition} \label{prop:bad:G23}
Let $G$ be a {\rm (}locally{\rm )} precompact abelian group such that 
$G[2]$ or $G[3]$ is sequentially open in $G$. Then $G$ admits no 
non-trivial quasi-convex null sequences. In particular, if $2G$ or $3G$ 
admits no non-trivial null sequences, then $G[2]$ or $G[3]$ is 
sequentially open in $G$, and hence $G$ admits no non-trivial 
quasi-convex null sequences.
\end{proposition}

\begin{proof} 
Suppose that $G[p]$ is sequentially open in $G$, where $p = 2$ or $p = 3$, 
and let $\{x_n\} \subseteq G$ be a~quasi-convex null sequence. Then $x_n 
\in G[p]$ for all but finitely many $n \in \mathbb{N}$.  Let $\{x_{n_k}\}$ 
denote the subsequence of $\{x_n\}$ consisting of all members that belong 
to~$G[p]$. By Corollary~\ref{cor:prel:subgrpA}, $\{x_{n_k}\}$ is 
a~quasi-convex null sequence in $G[p]$. Therefore, by 
Lemma~\ref{lemma:bad:Z23}, $\{x_{n_k}\}$~is trivial, and hence $\{x_n\}$ 
is trivial.

In order to show the second statement, we observe that $pG$ is a 
continuous homomorphic image of $G/G[p]$. Thus, if $pG$ has no non-trivial 
null sequences, then $G/G[p]$ has no non-trivial null sequences either, 
and therefore $G[p]$ is sequentially open in $G$.
\end{proof}

\begin{lemma} \label{lemma:bad:lift}
Let $\mathcal{P}$ be a topological property that is an inverse invariant 
of open perfect maps, $p$~a~prime number, and $D$ a non-trivial precompact 
abelian group of exponent $p$ with property $\mathcal{P}$. Then there 
exists a~minimal abelian group of exponent $p^2$ with property 
$\mathcal{P}$ such that $pG\cong D$.
\end{lemma}

\begin{proof}
Since $D$ has exponent $p$, so does its completion $\widetilde D$, and 
thus $\widetilde D \cong \mathbb{Z}_p^\lambda$ for some cardinal $\lambda$ 
(cf.~\cite[4.2.2]{DikProSto}). Put $K := \mathbb{Z}_{p^2}^\lambda$, and 
let $f\colon K \rightarrow pK \cong \widetilde D$ denote the continuous 
homomorphism defined by $f(x) = px$. Since $K$ is compact, the map $f$ is 
open and perfect. Consequently, $G := f^{-1}(D)$ has 
property~$\mathcal{P}$, the 
exponent of $G$ is $p^2$ (because $pG = f(G)= D$ is non-trivial), and $G$ 
is dense in~$K$. In particular, $K = \widetilde G$, and thus 
$\operatorname{soc}(\widetilde G) = K[p] = \ker f \subseteq G$. Therefore, 
by Corollary~\ref{cor:prel:st}(a), $G$ is minimal.
\end{proof}

\begin{proof}[Proof of Theorem~\ref{thm:main:ABSbad}.]
(a) Let $\kappa$ be an infinite cardinal, and let $D$ denote the direct 
sum $\mathbb{Z}_p^{(\kappa)}$ equipped with the Bohr-topology. By 
Lemma~\ref{lemma:bad:lift} (with $\mathcal{P}$ the trivial property), 
there exists a minimal abelian group $G$ of exponent $p^2$ such that $pG 
\cong D$. By a well-known theorem of Flor, the Bohr-topology of a discrete 
abelian group admits no non-trivial convergent sequences 
(cf.~\cite{Flor}). Hence, by Proposition~\ref{prop:bad:G23}, $G$ admits no 
non-trivial quasi-convex null sequences. This completes the proof, because 
$|pG| = |\mathbb{Z}_p^{(\kappa)}| = \kappa$.

(b) Let $\kappa$ be an infinite cardinal such that $\kappa^\omega=\kappa$. 
By a theorem of Dijkstra and van Mill (cf.~\cite[Theorem~3]{DijvMil}; 
see also \cite[5.8]{GaliMaca}),
the compact group $\mathbb{Z}_p^\kappa$ admits a subgroup $D$ such that:

\begin{mynumlist}

\item
$D$ contains no non-trivial convergent sequences; 

\item
$D$ is dense in the $G_\delta$-topology of $\mathbb{Z}_p^\kappa$;  and

\item
$|D|=\kappa$.

\end{mynumlist}

\smallskip

\noindent
It follows from property (2) that $D$ is pseudocompact 
(cf.~\cite[1.2]{ComfRoss2}). Since pseudocompactness is an inverse 
invariant of open perfect maps (cf.~\cite[3.10.H]{Engel6}), by 
Lemma~\ref{lemma:bad:lift}, there exists a pseudocompact minimal group $G$ 
of exponent $p^2$ such that $pG \cong D$.
Hence, by Proposition~\ref{prop:bad:G23}, $G$ admits no non-trivial
quasi-convex null sequences. This completes the proof, because
$|pG|=|D|=\kappa$ by property (3).
\end{proof}

\begin{proof}[Proof of Theorem~\ref{thm:main:MACHbad}.] 
Van Douwen showed that under MA, there exists an infinite countably 
compact abelian group $D$ of exponent $2$ that admits no non-trivial 
convergent sequences (cf.~\cite[8.1]{vDouw}), and he also observed that 
under CH, a~construction of Hajnal and Juh\'asz yields a group $D$ with 
the same properties (cf.~\cite{HajJuhCC} and~\cite[9]{vDouw}). Since 
countable compactness is an inverse invariant of perfect maps 
(cf.~\cite[3.10.10]{Engel6}), by Lemma~\ref{lemma:bad:lift}, there exists 
a countably compact minimal group $G$ of exponent $4$ such that $2G \cong 
D $. Hence, by Proposition~\ref{prop:bad:G23}, $G$ admits no non-trivial 
quasi-convex null sequences.
\end{proof}

\section{Bounded precompact groups without non-trivial quasi-convex null 
sequences}

\label{sect:BP}

\begin{samepage}
\begin{Ltheorem*}[\ref{thm:main:BP}]
For every bounded precompact abelian group $G$, the following statements 
are equivalent:

\begin{myromanlist}

\item
$G$ admits no non-trivial quasi-convex null sequences;

\item
one of the subgroups  $G[2]  = \{g \in G \mid 
2g =  0\}$ or $G[3]  = \{g  \in  G \mid  3g  =  0\}$ is sequentially open 
in~$G$.

\end{myromanlist}

\end{Ltheorem*}
\end{samepage}

In this section, we prove Theorem~\ref{thm:main:BP}. 
Recall that a set $\{f_1,\ldots,f_n\}$ 
of non-zero elements in an abelian group $G$ is {\itshape independent} if 
whenever $\sum\limits_{i=1}^n l_i f_i=0$ for some $l_i\in 
\mathbb{Z}$, then $l_i f_i=0$ for every $i$, or equivalently, if
$\langle f_1,\ldots,f_n \rangle = \langle f_1 \rangle \oplus \cdots
\oplus \langle f_n \rangle$. In what follows, $o(g)$ denotes the order 
of an element $g$ in a group.

\begin{lemma} \label{lemma:BP:indp}
Let $\{f_1,\ldots,f_n\}$ be an independent subset of a (locally)  
precompact abelian group~$G$. If $4 \leq o(f_i) < \infty$ for every $1\leq 
i \leq n$, then $X:=\{0\} \cup \{\pm f_1,\ldots, \pm f_n\}$ is 
quasi-convex in $G$.
\end{lemma}

\begin{proof}
Set $m_k:=o(f_k)$ for $1\leq k \leq n$,  $m_k=4$ for $k >n$, and
$P:=\prod\limits_{k=1}^\infty \mathbb{Z}_{m_k}$. For every $k\geq 1$,
put $e_k:=(0,\ldots,0,1,0,\ldots)$, with $1$ at the $k$-th coordinate and 
zero elsewhere. The authors showed in \cite[5.5]{DikGL2} that
$S:=\{0\}\cup\{\pm e_k \mid k\geq 1\}$ is quasi-convex in $P$.

Put $F:=\langle f_1,\ldots,f_n \rangle = \langle f_1 \rangle \oplus 
\cdots\oplus \langle f_n \rangle$.  Clearly, $F$ is finite. Consequently, 
the homomorphism  $\varphi \colon F \rightarrow P$ defined by 
$\varphi(f_i)=e_i$ for  $1 \leq i \leq n$ is an embedding of topological 
groups. Therefore, by Lemma~\ref{lemma:prel:qc-hom},
$X=\varphi^{-1}(S)$ is quasi-convex in $F$. Hence,  by 
Corollary~\ref{cor:prel:subgrpB} and  Lemma~\ref{lemma:prel:subgr}(b), $X$ 
is quasi-convex in $G$.
\end{proof}

\begin{proposition} \label{prop:BP:FEomega}
Let $E$ and $F$ be finite abelian groups, and suppose that 
$\exp E \geq 4$. Then every dense subgroup $A \leq F \times E^\omega$ 
contains a non-trivial null sequence that
is quasi-convex both in $A$ and in $F \times E^\omega$.
\end{proposition}

\begin{proof}
For every positive integer $n$, 
let $\pi_n\colon F \times E^\omega \rightarrow F \times E^n$ denote 
the canonical projection of the first $n+1$ coordinates. 
Pick $y\in E$ such that $o(y)=\exp E$.  Since $A$ is 
dense in $F \times E^\omega$, one has $\pi_n(A)=F\times E^n$. Thus, 
for every $n$, we may pick $x_n \in A$ such that
$\pi_n(x_n)=(0,\ldots,0,y)$. We claim that $\{x_n\}$ is a 
quasi-convex null sequence in $A$ and $F \times E^\omega$.

{\itshape Step 1.}
We show by induction on $n$ that  the set 
$\{\pi_n(x_1),\ldots,\pi_n(x_n)\}$ is independent in 
$F\times E^n$.  For $n=1$, the statement is trivial, 
because $\pi_1(x_1)$ is 
non-zero. Assume now that the statement holds for~$n$, and suppose that
$\sum\limits_{i=1}^{n+1} l_i \pi_{n+1}(x_i)=0$ for 
$l_i \in \mathbb{Z}$. Then
$\sum\limits_{i=1}^{n+1} l_i \pi_{n}(x_i)=0$, and thus
$\sum\limits_{i=1}^{n} l_i \pi_{n}(x_i)=0$,
because $\pi_n(x_{n+1})=0$. By the 
inductive hypothesis, it follows that
$l_i \pi_{n}(x_i)=0$ for 
$1\leq i \leq n$.  The $(i+1)$\mbox{-}th coordinate of $\pi_{n}(x_i)$ is 
$y$, 
and so $o(y) \mid l_i$ for $1\leq i \leq n$. Therefore,
$l_i x_i=0$ for $1\leq i \leq n$, because $o(y)=\exp E$. 
Hence, $l_{n+1} \pi_{n+1}(x_{n+1}) = 
- \sum\limits_{i=1}^n l_i \pi_{n+1}(x_i) = 0$,
as required.

{\itshape Step 2.}
Put $S :=  \{0\}\cup\{\pm x_n \mid n  \in \mathbb{N}\}$. By 
Lemma~\ref{lemma:BP:indp},
$\pi_n(S) = \{0\} \cup \{\pm \pi_{n}(x_1),\ldots, \pm 
\pi_{n}(x_n)\}$ is quasi-convex in $F \times E^n$, because
$o(\pi_n(x_i))=o(y) = \exp E \geq 4$. Thus, by  
Lemma~\ref{lemma:prel:qc-hom},
\begin{align} 
Q_{F \times E^\omega}(S) 
\subseteq \pi_n^{-1}(\pi_n(S)) = S + \ker \pi_n
\end{align}
for every $n$. Therefore,
\begin{align} \label{dense:eq:QS'}
Q_{F\times E^\omega}(S) \subseteq \bigcap\limits_{n=1}^\infty (S + \ker 
\pi_n) = 
\operatorname{cl}_{F \times E^\omega} S,
\end{align}
because $\{\ker \pi_n\}_{n=1}^\infty$ is a base for the topology of 
$F \times E^\omega$ at zero.
Since $x_k \in \ker \pi_n$ for every $k>n$, it 
follows that $\{x_n\}$ is a null sequence, and 
$S$ is closed in $F \times E^\omega$. Hence, by \eqref{dense:eq:QS'}, 
$S$ is quasi-convex in~$F \times E^\omega$, as desired. By 
Corollary~\ref{cor:prel:subgrpA}, this implies that $\{x_n\}$ is 
quasi-convex in $A$ as well.
\end{proof}

\begin{corollary} \label{cor:BP:AK}
Let $K$ be an infinite bounded compact metrizable abelian group that 
contains no open compact subgroup of the form $\mathbb{Z}_2^\omega$
and $\mathbb{Z}_3^\omega$. Then every dense subgroup $A$ of $K$ contains 
a~non-trivial null sequence that is quasi-convex in both  $A$ and $K$.

\end{corollary}

\begin{proof}
Since $K$ is bounded compact abelian, it is topologically isomorphic to a 
product of finite cyclic groups (cf.~\cite[4.2.2]{DikProSto}). The 
number of the factors is countably infinite, because $K$ is metrizable and 
infinite, and the number of non-isomorphic factors is finite, as $K$ is 
bounded. Thus,
$K \cong F \times \mathbb{Z}_{m_1}^\omega \times \cdots\times
 \mathbb{Z}_{m_l}^\omega$, where $F$ is a finite abelian group
and $m_1,\ldots,m_l$ are distinct integers. Consequently, for 
$E:=\mathbb{Z}_{m_1} \times \cdots \times \mathbb{Z}_{m_l}$, one has
$K \cong F \times E^\omega$. Since $K$ contains no open compact subgroup 
of the form $\mathbb{Z}_2^\omega$ and $\mathbb{Z}_3^\omega$,
clearly $E \neq \mathbb{Z}_2$ and $E \neq \mathbb{Z}_3$. Therefore,
$\exp(E) \geq 4$, and hence the statement follows by 
Proposition~\ref{prop:BP:FEomega}.
\end{proof}

\begin{corollary} \label{cor:BP:A23}
Let $A$ be a bounded precompact metrizable abelian group. If the subgroups 
$A[2]$ and $A[3]$ are not open in $A$, then $A$ contains a non-trivial 
null sequence that is quasi-convex in both $A$ and the completion 
$\widetilde A$ of $A$.
\end{corollary}

\begin{proof}
Put $K:=\widetilde A$. Since $A$ is metrizable and bounded, so is $K$. One 
has $A[p]=K[p] \cap A$, and thus $K[2]$ and $K[3]$ are not open in $K$. In 
particular, $K$ is infinite, and it contains no open compact subgroup of 
the form $\mathbb{Z}_2^\omega$ and $\mathbb{Z}_3^\omega$. Hence, the 
statement follows from Corollary~\ref{cor:BP:AK}.
\end{proof}

\begin{lemma}\label{lemma:BP:seqmet} 
If $A$ is a bounded precompact  abelian group  generated by a null sequence 
$\{w_n\}$, then $A$ is metrizable.
\end{lemma}

\begin{proof} 
Let $m$ denote the exponent of $A$. Clearly, every character $\chi \in 
\widehat{A}$ is completely determined by the values $\chi(w_n)$ taken at 
the generators of $A$. Since $\{w_n\}$ is a null sequence and $\chi$ is 
continuous,  $\chi(w_n)\to 0$ in $\mathbb{T}$, and $\chi(x_n)$  belongs 
to the cyclic subgroup of order $m$ in $\mathbb{T}$, because $mA=0$.
Consequently, $\chi(x_n)=0$ for all but finitely many indices $n$. 
Therefore, $\widehat{A}$ is countable. Hence, $A$ is metrizable
(cf.~\cite[2.12]{HernMaca}).
\end{proof}

\begin{proof}[Proof of Theorem~\ref{thm:main:BP}.] 
(ii) $\Rightarrow$ (i): This implication  holds even without the 
assumption that $G$ is bounded, and has already been shown in 
Proposition~\ref{prop:bad:G23}.

(i) $\Rightarrow$ (ii): Suppose that neither $G[2]$ nor 
$G[3]$ is sequentially open in $G$. Then there are sequences $\{y_n\}$ and 
$\{z_n\}$ in $G$ that witness that $G[2]$ and $G[3]$ are not sequentially 
open. In other words, $y_n \rightarrow y_0 \in G[2]$, but 
$y_n \notin G[2]$ for infinitely many indices $n$, and 
$z_n \rightarrow z_0 \in G[3]$, but $z_n \notin G[3]$ for infinitely many 
indices $n$. By replacing $\{y_n\}$ with $\{y_n  - y_0\}$ and $\{z_n\}$ 
with $\{z_n - z_0\}$, we may assume that $y_n \rightarrow 0$ and 
$z_n \rightarrow 0$. 

Let $\{w_n\}$ denote the alternating sequence $y_1,z_1,y_2,z_2,\ldots$. 
Clearly, $\{w_n\}$ is a null sequence, and $w_n \notin G[2]$  for 
infinitely many indices $n$ and $w_n \notin G[3]$ for infinitely many 
indices $n$. Let $A$ denote the subgroup of $G$ generated by $\{w_n\}$.
Then $A$ is a bounded precompact abelian group, and by 
Lemma~\ref{lemma:BP:seqmet}, $A$ is metrizable. Although 
$w_n \rightarrow 0$, one has $w_n \notin A[2]=G[2]\cap A$ for infinitely 
many indices $n$ and $w_n \notin A[3]=G[3]\cap A$ for infinitely many
indices $n$. Thus, the subgroups $A[2]$ and $A[3]$ are not open in $A$.
Therefore, by Corollary~\ref{cor:BP:A23}, there is a non-trivial null 
sequence $\{x_n\}\subseteq A$ such that $\{x_n\}$ is quasi-convex both in 
$A$ and the completion $\widetilde A$ of $A$.

Since the completion $\widetilde A$ is a closed subgroup of the completion 
$\widetilde G$ of $G$, by Corollary~\ref{cor:prel:subgrpB}, 
$\{x_n\}$ is also quasi-convex in $\widetilde G$. Hence, by 
Corollary~\ref{cor:prel:subgrpA}, $\{x_n\}$ is quasi-convex in $G$, 
because $\{x_n\} \subseteq A \subseteq G$.
\end{proof}

\section{Sequences of the form 
$\boldsymbol{\{\frac 1 {b_n} \}_{n=0}^\infty}$ in $\boldsymbol{\mathbb{T}}$}

\label{sect:T}

\begin{Ltheorem*}[\ref{thm:main:8T}] Let $\{q_n\}_{n=0}^\infty$ be a sequence of positive integers, and put
$b_n =  q_0\cdots q_n$ \ for every  $n \in \mathbb{N}$.
If $q_n \geq  8$ for every $n \in \mathbb{N}$, then $\{\frac 1 {b_n} \}_{n=0}^\infty$ is a quasi-convex sequence  in $\mathbb{T}$.
\end{Ltheorem*}

In this section, we prove Theorem~\ref{thm:main:8T}.  Although in 
Theorem~\ref{thm:main:8T} itself we require $q_n \geq 8$, a number of 
intermediate statements remain true under no conditions at all or weaker 
conditions imposed upon the sequence $\{q_n\}$. Consequently, we consider 
$\{q_n\}_{n=0}^\infty$ (and thus $\{b_n\}_{n=0}^\infty$) a fixed sequence 
of positive integers and set 
$X := \{0\} \cup  \{\pm \frac 1 {b_n} \mid n \in \mathbb{N}\}$ throughout 
this section, but 
make no further  assumptions about their properties; instead, we impose 
conditions on the $\{q_n\}$ in each statement as needed.

The first step toward the proof of Theorem~\ref{thm:main:8T} is to 
establish a standard method for describing elements in $Q_\mathbb{T}(X)$.
As we have seen in 
\cite{DikGL1} and \cite{DikGL2}, finding a 
convenient way to represent elements of $\mathbb{T}$ is a useful tool in 
calculating quasi-convex hulls of sequences. Let $\{d_i\}_{i=0}^\infty$ be 
an increasing sequence of positive integers such that $d_i \mid d_{i+1}$ 
for every $i \in \mathbb{N}$. Then 
$z \in \mathbb{T}$ (which, as we stated 
in the Introduction, is 
identified with $(-\frac 1 2, \frac 1 2]$) can be expressed in the form
$z= \sum\limits_{i=0}^\infty \frac{c_i}{d_i}$, 
with $c_i$ integers such that 
$|c_i| \leq \frac{d_i}{2d_{i-1}}$. (We 
consider $d_{-1}  =  1$.) 
This representation, however, need not be unique: For  example, if 
$d_0  =  3$ and 
$d_1  =  6$, then $\frac 1 6$ can be 
expressed with 
$c_0  = 0$ and 
$c_1 = 1$, but also with 
$c_0  = 1$ and 
$c_1 = -1$. In order to eliminate this 
anomaly, we say that the representation of $z$ is {\itshape standard} if 
the following conditions are satisfied:

\begin{myromanlist}

\item
$c_i  \in  \mathbb{Z}$ and 
$|c_i|  \leq  \frac{d_i}{2d_{i-1}}$ for  all  $i  \in   \mathbb{N}$;

\item
$\left|z- \sum\limits_{i=0}^k \frac{c_i}{d_i}\right|
  \leq  \frac 1 {2d_k}$ for every 
$k  \in  \mathbb{N}$;

\item
if 
$\left|z- \sum\limits_{i=0}^k \frac{c_i}{d_i}\right|  =  \frac 1 {2d_k}$ for some $k$, then 
$\left|\frac{c_k}{d_k}\right| 
  <   
\left|z- \sum\limits_{i=0}^{k-1} \frac{c_i}{d_i}\right|$.

\vspace{6pt}

\end{myromanlist}

\noindent
(In the aforementioned example, $c_0   = 0$ and $c_1  = 1$ is a standard representation of $\frac 1 6$, but $c_0 = 1$ and $c_1  =  -1$ 
is not a standard one.)
The next theorem is in the vein of \cite[4.3]{DikGL1} and 
\cite[3.2]{DikGL2}, and plays an important role in the proof of 
Theorem~\ref{thm:main:8T}.

\begin{theorem} \label{thm:T:-101}
If $q_{k+1}  \geq  4$ whenever $q_k  =  7$, then
$Q_\mathbb{T}(X) \subseteq 
\left\{\sum\limits_{i=0}^\infty \dfrac{\varepsilon_i} {b_i} \mid
\varepsilon_i \in \{-1,0,1\}\right\}$. Furthermore, every 
$x\in Q_\mathbb{T}(X)$ admits a standard representation
$\sum\limits_{i=0}^\infty \dfrac{\varepsilon_i} {b_i}$ with
$\varepsilon_i \in \{-1,0,1\}$.
\end{theorem}

The proof of Theorem~\ref{thm:T:-101} requires two preparatory steps. The 
first one is an analogue  of \cite[2.1, 2.2, 2.4]{DikGL2}, for which we 
introduce our own rounding functions: For 
$x  \in \mathbb{R}$, we put
\begin{align*}
\lceil x \rceil & := \min \{n \in \mathbb{Z} \mid  x < n\}, & 
[x] & := \max \{n \in \mathbb{Z} \mid  n \leq x\},& 
\lfloor x \rfloor & := \max \{n \in \mathbb{Z} \mid  n < x\}. 
\end{align*}
We note that these are not the usual definitions of the floor and ceiling functions (as we use strict inequality in both).

\begin{lemma} \label{lemma:T:repz}
Let $z = \sum\limits_{i=0}^\infty \frac{c_i}{d_i}  \in  \mathbb{T}$ be a standard representation.

\begin{myalphlist}

\item
If $mz  \in \mathbb{T}_+$ for all $m  =  1,\ldots,\lceil \frac{d_0}{6} \rceil$, then $c_0 \in \{-1,0,1\}$.

\item
If $mz \in \mathbb{T}_+$ for all  $m =   1,\ldots,[\frac{d_0}{4}]$ and $d_0  \neq  7$, then $c_0  \in  \{-1,0,1\}$.

\item
If $mz \in \mathbb{T}_+$ for all $m  =  1,\ldots,[\frac{d_0}{4}]$ and for $m =  d_0-1$, then $c_0 \in \{-1,0,1\}$.

\end{myalphlist}
\end{lemma}

\begin{proof} (a) Put $l  =  \lceil\frac{d_0}{6}\rceil$. Since $mz \in \mathbb{T}_+$ for  all
$m  =  1,\ldots,l$, one has $z \in \{1,\ldots,l\}^\triangleleft 
 =  \mathbb{T}_l$, and thus $|z| \leq  \frac 1 {4l}  <    \frac{3}{2d_0}$. Therefore,
\begin{align}
\left| \frac{c_0}{d_0} \right|
\leq |z| +  \left| z - \frac{c_0}{d_0}\right| < 
\frac{3}{2d_0}+ \frac 1 {2{d_0}} =\frac{2}{d_0}.
\end{align}
Hence, $|c_0|  <  2$, as desired.

\pagebreak[3]

(b) For $d_0  = 2$ and 
$d_0 = 3$,  the conclusion is trivial. 
If $d_0  \neq  6,7$, then one has 
$\lceil \frac{d_0} 6 \rceil  \leq 
[\frac{d_0} 4]$, and the statement follows from part (a).
Suppose that $d_0 = 6$. Then it is 
given that  
$z \in  \mathbb{T}_+$, and thus 
$|z|  \leq   \frac 1 4$. If 
$|z-\frac{c_0}{d_0}|  <   \frac{1}{2d_0}$, 
then
\begin{align}
\left|\frac {c_0}{d_0}\right| \leq |z| + \left|z-\frac {c_0}{d_0}\right| 
< \frac 1 4 + \frac 1 {2d_0} = \frac{2}{d_0},
\end{align}
and thus $|c_0|   <  2$. 
If $|z-\frac {c_0}{d_0}|  
=   \frac{1}{2d_0}$, then 
$|\frac {c_0}{d_0}|   < 
|z|  \leq  \frac 1 4$, because 
the  representation of $z$ is standard. Hence, 
$|c_0|  <  2$, as desired.

(c) If $d_0\neq  7$, 
then the statement follows from (b), and so we 
may suppose that $d_0 = 7$. Then it 
is given that  
$z,6z \in \mathbb{T}_+$,  
which means that
\begin{align}
z \in \{1,6\}^\triangleleft = \mathbb{T}_6 \cup 
(-\tfrac 1 6+\mathbb{T}_6) \cup (\tfrac 1 6+\mathbb{T}_6).
\end{align}
Thus, $|z| \leq \tfrac 5 {24}$. 
Therefore,
\begin{align}
\left| \dfrac{c_0}{d_0} \right| \leq
|z| + \left|z - \dfrac{c_0}{d_0} \right|  \leq
\dfrac{5}{24} + \dfrac{1}{14} = \dfrac{47}{168} < \dfrac{48}{168} = 
\dfrac{2}{d_0}.
\end{align}
Hence, $|c_0|  <  2$.
\end{proof}

The second preparatory step to precede the proof of  
Theorem~\ref{thm:T:-101} involves finding characters in 
$X^\triangleright$. Let 
$\eta_0\colon \mathbb{T} \rightarrow \mathbb{T}$ 
denote the  identity homomorphism, and for 
$k \geq 1$, set 
$\eta_k  :=  b_{k-1} \eta_0$.

\begin{lemma} \label{lemma:T:etak}
If $q_k \geq 4$, then 
$m\eta_k \in X^\triangleright$ for
$m = 
1,\ldots,[\frac{q_k}{4}]$. If in addition 
$q_{k+1}\geq 4$, then 
$m\eta_k \in X^\triangleright$  also for 
$m = q_k -1$.  
\end{lemma}

\begin{proof}
Fix $n \in \mathbb{N}$. 
If $n < k$, then 
$\eta_k(\frac 1 {b_n}) = 
\frac{b_{k-1}}{b_n} \equiv_1 0$, and thus
$m\eta_k(\frac 1 {b_n}) \in
 \mathbb{T}_+$. 

If $n= k$, then
$\eta_k(\frac 1 {b_n})  = 
\frac{b_{k-1}}{b_k} =  \frac 1 {q_k}$, 
and so
$m\eta_k(\frac 1 {b_n})  = 
\frac m {q_k} \in \mathbb{T}_+$ for
$m =  1,\ldots,[\frac{q_k}{4}]$ and  
$m  = q_k -1$.

If  $n  >  k$, then 
$\eta_k(\frac 1 {b_n})  = 
\frac{b_{k-1}}{b_n}  =  
\frac{1}{q_k q_{k+1}\cdots q_n}$ and
$m\eta_k(\frac 1 {b_n})  =
\frac{m}{q_k q_{k+1}\cdots q_n}$. 
Consequently,
$m\eta_k(\frac 1 {b_n})  \in \mathbb{T}_+$ for 
$m  = 1,\ldots,[\frac{q_k}{4}]$. If 
$q_{k+1}\geq  4$, then
$\frac{q_k-1}{q_k q_{k+1}\cdots q_n} \leq  
 \frac {1}{q_{k+1}} \leq \frac 1 4$, 
and hence 
$(q_k-1)\eta_k(\frac 1 {b_n}) \in
 \mathbb{T}_+$.
\end{proof}

We are now ready to prove Theorem~\ref{thm:T:-101}.

\begin{proof}[Proof of Theorem~\ref{thm:T:-101}.]
Let $x \in Q_\mathbb{T}(X)$, and let
$x  = \sum\limits_{i=0}^\infty 
\frac{c_i}{b_i}$ be a standard representation of $x$. 
Fix $k \in \mathbb{N}$, and put 
$d_i  :=  \frac{b_{k+i}}{b_{k-1}} =  q_k \cdots q_{k+i}$ 
for every $i  \in \mathbb{N}$. (As usual, we consider $b_{-1}  = 1$.)
Then $z := \eta_k(x) =  b_{k-1}x  \equiv_1 \sum\limits_{i=0}^\infty 
\frac{c_{k+i}}{d_i}$, and it is  a standard representation of $z$, 
because $x  = \sum\limits_{i=0}^\infty \frac{c_i}{b_i}$ is a standard
representation of $x$.  (Indeed, 
$|c_{k+i}| \leq \frac{b_{k+i}}{2 b_{k+i-1}} = \frac{d_i}{2 d_{i-1}}$, 
while conditions (ii) and (iii) follows by observing that
$\left| z - \sum\limits_{i=0}^m \frac{c_{k+i}}{d_i}\right|=
b_{k-1}\left|x -\sum\limits_{i=0}^{k+m} \frac{c_i}{b_i}\right|$ for every 
$m\in \mathbb{N}$.)
Furthermore, if  $m\eta_k  \in  X^\triangleright$, 
then $mz  = m\eta_k(x) \in \mathbb{T}_+$.

If $q_k <  4$, then $|c_k|  \leq   \frac{q_k}{2} <  2$, and there is 
nothing to prove.  So, we may assume that 
$q_k  \geq  4$. Thus, by 
Lemma~\ref{lemma:T:etak}, 
$m\eta_k \in X^\triangleright$ for 
$m  =  1,\ldots,[\frac{q_k}{4}]$. 
Consequently, $m z  \in \mathbb{T}_+$ for 
$m  =  1,\ldots,[\frac{q_k}{4}]$. 
Therefore,
if $q_k  \neq  7$, then by 
Lemma~\ref{lemma:T:repz}(b), 
the first coefficient of $z$, that is $c_k$, satisfies
$c_k \in  \{-1,0,1\}$.
If $q_k =  7$,  then by our 
assumption, $q_{k+1} \geq  4$, and 
by  Lemma~\ref{lemma:T:etak}, 
$(q_k-1)\eta_k \in 
X^\triangleright$.
Consequently, by Lemma~\ref{lemma:T:repz}(c), the first coefficient of 
$z$, that is $c_k$, satisfies
$c_k \in \{-1,0,1\}$.
\end{proof}

The next lemma is somewhat similar to \cite[3.3]{DikGL2}, both in its 
content and its role in the proof of Theorem~\ref{thm:main:8T}.

\begin{lemma}  \label{lemma:T:k1k2}
Let $k_1,k_2 \in \mathbb{N}$ be such that 
$k_1 < k_2$. 
Then $[\frac{q_{k_1}} 4]\eta_{k_1}
 \pm 
\lfloor \frac{q_{k_2}} 4\rfloor
\eta_{k_2}  \in  
X^\triangleright$.
\end{lemma}

\begin{proof}
Let $n \in\mathbb{N}$.  
If $n  <  k_2$, then 
$\eta_{k_2}(\frac 1 {b_n})  =  
\frac{b_{k_2-1}}{b_n} \equiv_1  0$, 
and thus 
\begin{align}
([\frac{q_{k_1}} 4]\eta_{k_1}\pm \lfloor \frac{q_{k_2}} 
4\rfloor\eta_{k_2}) (\frac 1 {b_n}) 
& \equiv_1 
[\frac{q_{k_1}} 4]\eta_{k_1} (\frac 1 {b_n}) \in \mathbb{T}_+,
\end{align}
because 
$[\frac{q_{k_1}} 4]  \eta_{k_1} \in 
X^\triangleright$
by Lemma~\ref{lemma:T:etak}. Suppose now that 
$k_2  \leq  n$. Then
\begin{align}
|[\frac{q_{k_1}} 4] \eta_{k_1}(\frac 1 {b_n})| & \leq
\frac{[\frac{q_{k_1}} 4] }{q_{k_1} q_{k_2}} \leq 
\frac 1 {4q_{k_2}},  \mbox{ and} \\
|\lfloor \frac{q_{k_2}} 4\rfloor \eta_{k_2}(\frac 1 {b_n})| & \leq 
\frac{ \lfloor \frac{q_{k_2}} 4\rfloor}{q_{k_2}} \leq 
\frac 1 4 - \frac{1}{4q_{k_2}}.
\end{align}
Therefore, 
\begin{align}
|([\frac{q_{k_1}} 4]\eta_{k_1}\pm \lfloor \frac{q_{k_2}} 
4\rfloor\eta_{k_2})
(\frac 1 {b_n})| & \leq 
|[\frac{q_{k_1}} 4] \eta_{k_1}(\frac 1 {b_n})| +
|\lfloor \frac{q_{k_2}} 4\rfloor \eta_{k_2}(\frac 1 {b_n})| \leq 
\frac 1 4,
\end{align}
and hence 
$([\frac{q_{k_1}} 4]\eta_{k_1} \pm 
\lfloor \frac{q_{k_2}} 
4\rfloor\eta_{k_2}) (\frac 1 {b_n}) 
\in \mathbb{T}_+$.
This shows that 
$[\frac{q_{k_1}}4]\eta_{k_1}
\pm \lfloor\frac{q_{k_2}}4\rfloor\eta_{k_2}
\in X^\triangleright$.
\end{proof}

We introduce further notations to facilitate calculations  in 
$\mathbb{T}$. Let
$z  = 
\sum\limits_{i=0}^\infty \frac{c_i}{d_i}$ 
be a~standard 
representation of $z \in \mathbb{T}$ with 
respect to a sequence
$\{d_i\}_{i=0}^\infty$ such that $d_i \mid d_{i+1}$
for every $i \in  \mathbb{N}$. We put
$\Lambda(z)  := 
\{i \in \mathbb{N} \mid 
c_i  \neq  0\}$,
$\underline{q}(z)   :=  
\min \{ \frac{d_{i+1}} {d_i} \mid i \in  
\Lambda(z)\}$, and 
$S(x)  :=  \frac{1}{\underline q(x) - 1}$.
(As usual, we consider 
$d_{-1}  =  1$.)

\begin{lemma} \label{lemma:T:qLambda}
Let $x  = \sum\limits_{i=0}^\infty \frac{\varepsilon_i}{b_i}$ be a
standard representation of $x \in \mathbb{T}$, with  $\varepsilon_i \in \{-1,0,1\}$.  
Then, for every $k \in \Lambda(x)$, one has $\frac 1 {q_k}(1-S(x))  \leq  |\eta_k(x)|  \leq  \frac 1 {q_k} (1+S(x))$.
\end{lemma}

\begin{proof}
One has $\eta_k(x) = b_{k-1} x \equiv_1 \sum\limits_{i=0}^\infty 
\frac{\varepsilon_{k+i}}{q_k\cdots q_{k+i}}$,
and thus one obtains (modulo $1$)
\begin{align}
\left| \eta_k(x) - \frac{\varepsilon_{k}}{q_{k}} \right| \leq
\frac{1}{q_k} \sum\limits_{i=1}^\infty \frac{1}{q_{k+1}\cdots 
q_{k+i}} \leq
\frac 1 {q_k} \sum\limits_{i=1}^\infty
\frac{1}{(\underline q(x))^{i}} = 
\frac {S(x)} {q_k}.
\end{align}
Therefore, 
$\frac 1 {q_k}(1-S(x)) 
 \leq  |\eta_k(x)| 
 \leq  \frac 1 {q_k} (1+S(x))$,
as required.
\end{proof}

\begin{lemma} \label{lemma:T:estimate}
Suppose that $q_{k+1}  \geq  4$ 
whenever $q_k  = 7$, and let $x \in Q_\mathbb{T}(X)$. If
$k_1,k_2 \in \Lambda(x)$ are such that $k_1 < k_2$, then
\begin{align}
\left(\frac{[\frac{q_{k_1}} 4]}{q_{k_1}} + 
\frac{\lfloor \frac{q_{k_2}} 4\rfloor }{q_{k_2}}\right)(1-S(x)) \leq
\frac 1 4.
\end{align}
\end{lemma}

\begin{proof}
Put $\chi  := \varepsilon_{k_1}[\frac{q_{k_1}} 4]\eta_{k_1} + \varepsilon_{k_2} \lfloor \frac{q_{k_2}} 4\rfloor \eta_{k_2}$. 
By Lemma~\ref{lemma:T:k1k2}, 
$\chi \in X^\triangleright$, and so $\chi(x) \in \mathbb{T}_+$. The conditions imposed upon 
$\{q_n\}_{n=0}^\infty$ guarantee that Theorem~\ref{thm:T:-101} is applicable, and thus $x$ can be expressed in standard form as
$x  = \sum\limits_{i=0}^\infty \frac{\varepsilon_i}{b_i}$, with
$\varepsilon_i  \in \{-1,0,1\}$. In what follows, we use Lemma~\ref{lemma:T:qLambda} to estimate $\chi(x)$ from below.

If $\underline q(x) =  2$ or $\underline q(x)  =  3$, then the  statement is  trivial, and so we may  assume that $\underline q(x) \geq  4$. Then, by Lemma~\ref{lemma:T:qLambda}, $|\eta_{k_j}(x)| \leq \frac{4}{3q_{k_j}}$ ($j = 1,2$), and thus
\begin{align}
\lfloor \frac{q_{k_j}}{4} \rfloor |\eta_{k_j}(x)| \leq
[\frac{q_{k_j}}{4}] |\eta_{k_j}(x)| \leq
\frac {4[\frac{q_{k_j}}{4}]}{3 q_{k_j}} \leq \frac 1 3.
\end{align}
Therefore, 
\begin{align}
|\chi(x)|  \leq [\frac{q_{k_1}}{4} ]  |\eta_{k_1}(x)| + 
\lfloor \frac{q_{k_2}}{4} \rfloor  |\eta_{k_2}(x)| \leq \frac 2 3.
\end{align}
This implies that we can perform the  remaining calculations in 
$[-\frac 2 3,\frac 2 3] \subseteq\mathbb{R}$, and 
$\chi(x) \in \mathbb{T}_+$ if and only if 
$-\frac  1 4  \leq \chi(x) \leq  \frac  1 4$.
Since the first term of $\varepsilon_{k_j}\eta_{k_j}(x)$ 
is~$\frac{1}{q_{k_j}}$, one has $0 \leq \varepsilon_{k_j}\eta_{k_j}(x)$, 
because $|\varepsilon_{k_j}\eta_{k_j}(x) - \frac{1}{q_{k_j}}|\leq 
\frac{1}{2 q_{k_j}}$.
Consequently, by
Lemma~\ref{lemma:T:qLambda}, $\frac 1 {q_{k_j}}(1-S(x)) \leq 
\varepsilon_{k_j}\eta_{k_j}(x)$, and hence
\begin{align}
\left(\frac{[\frac{q_{k_1}} 4]}{q_{k_1}} + 
\frac{\lfloor \frac{q_{k_2}} 4\rfloor }{q_{k_2}}\right) (1-S(x)) \leq
(\varepsilon_{k_1}[\frac{q_{k_1}} 4]\eta_{k_1} + 
\varepsilon_{k_2} \lfloor \frac{q_{k_2}} 4\rfloor \eta_{k_2})(x) =
\chi(x) \leq 
\frac 1 4   ,
\end{align}
as desired.
\end{proof}

We are now ready to prove Theorem~\ref{thm:main:8T}.

\begin{proof}[Proof of Theorem~\ref{thm:main:8T}.] Let $x \in Q_\mathbb{T}(X)$, and assume that $x \notin X$.
Then $|\Lambda(x)|  >1 $, and so we may pick $k_1,k_2 \in \Lambda(x)$ such that $k_1  <  k_2$. As $q_n \geq 8$ for every $n\in\mathbb{N}$, one has
$1-S(x) \geq \frac{6}{7}$, and thus, by Lemma~\ref{lemma:T:estimate},
\begin{align}
\frac{[\frac{q_{k_1}} 4]}{q_{k_1}} + 
\frac{\lfloor \frac{q_{k_2}} 4\rfloor }{q_{k_2}} \leq \frac{7}{24}.
\end{align}
This inequality, however, does not hold with any $q_{k_j} \geq  8$. Hence, $x  \in  X$, as desired. 
\end{proof}

\begin{example} \label{ex:T:socT}
Let $\{p_n\}_{n=0}$ be an enumeration of all primes greater than~$8$, and  put $b_n  = p_0 \cdots p_n$ for every 
$n  \in  \mathbb{N}$. By Theorem~\ref{thm:main:8T}, $\{\frac 1 {b_n}\}_{n=0}^\infty$ is
quasi-convex in $\mathbb{T}$, and since each $b_n$ is square-free, $\{\frac 1 {b_n}\}_{n=0}^\infty \subseteq   \operatorname{soc}(\mathbb{T})$.
\end{example}

Using the next lemma, one can lift Theorem~\ref{thm:main:8T} into $\mathbb{R}$. Recall that in this note,  $\pi\colon \mathbb{R} \rightarrow \mathbb{T}$  denotes  the  canonical projection.

\begin{flemma}[{\cite[2.4]{DikGL1}}] \label{prel:lemma:pi-Y}
Let $Y \subseteq  \mathbb{R}$  be a compact subset. If there is $\alpha  \neq  0$ such that $\alpha Y   \subseteq (-\frac 1 2, \frac 1 2)$ and $\pi(\alpha Y)$ is quasi-convex in $\mathbb{T}$, then $Y$ is quasi-convex in $\mathbb{R}$. 
\end{flemma}

\begin{corollary}
Let $\{x_n\}_{n=0}^\infty \subseteq \mathbb{R}$ be a null sequence in $\mathbb{R}$ such that $q_n := \frac{x_{n-1}}{x_{n}}$ are  integers  and  $q_n   \geq  8$ for every $n \in  \mathbb{N}\backslash\{0\}$.  Then $\{x_n\}_{n=0}^\infty$ is quasi-convex in $\mathbb{R}$.
\end{corollary}

\begin{proof} 
Put $\alpha  = \frac{1}{8x_0}$, $q_0 = 8$, and $b_n  = q_0\cdots q_n$. 
Then  $\alpha x_n =   \frac{1}{b_n}$, and thus, by  
Theorem~\ref{thm:main:8T}, the sequence $\{\pi(\alpha x_n)\}_{n=0}^\infty$ 
is quasi-convex in $\mathbb{T}$.  Since 
$\{\alpha x_n\}_{n=0}^\infty \subseteq  [-\frac 1 8 ,\frac 1 8]$, by  
Lemma~\ref{prel:lemma:pi-Y}, the sequence $\{x_n\}_{n=0}^\infty$ is 
quasi-convex in $\mathbb{R}$, as required.
\end{proof}

\section{Compact-like abelian groups that admit no non-trivial 
quasi-convex  null sequences}

\label{sect:omega}

\begin{Ltheorem*}[\ref{thm:main:minimal}]
For every minimal abelian group $G$, the following statements are  equivalent:

\begin{myromanlist}
\item
$G$ admits no non-trivial quasi-convex null sequences;

\item
$G \cong P \times F$, where $P$ 
is a minimal bounded abelian $p$-group {\rm (}$p \leq 3${\rm )} admitting 
no non-trivial quasi-convex null sequences, and $F$ is a finite abelian 
group;

\item
one of the subgroups  $G[2]  = \{g \in G \mid 2g =  0\}$ or 
$G[3]  = \{g  \in  G \mid  3g  =  0\}$ is sequentially open in~$G$;

\item
$G$ contains a sequentially open compact subgroup of the form 
$\mathbb{Z}_2^\kappa$ or $\mathbb{Z}_3^\kappa$ \ for some cardinal $\kappa$.

\end{myromanlist}
\end{Ltheorem*}

\begin{Ltheorem*}[\ref{thm:main:omega}]
The following statements are equivalent for every abelian group $G$ that 
is  $\omega$-bounded  \underline{or} totally minimal:

\begin{myromanlist}

\item
$G$ admits no non-trivial quasi-convex null sequences;

\item
one of the subgroups $G[2]  = \{g \in G \mid 2g =  0\}$ or 
$G[3]  = \{g  \in  G \mid  3g  =  0\}$ is open in $G$;

\item
one of the subgroups $2G$ and $3G$ is finite.

{\parindent -25pt 
Furthermore, if $G$ is totally minimal, then these conditions are also 
equivalent to:
}

\item
$G  \cong \mathbb{Z}_2^\kappa  \times  F$ or 
$G \cong  \mathbb{Z}_3^\kappa  \times  F$, 
where $\kappa$ is some cardinal and $F$ is a finite abelian group.

\end{myromanlist}
\end{Ltheorem*}

In this section, we present the proofs of Theorems~\ref{thm:main:minimal} 
and~\ref{thm:main:omega}. Since the latter relies on the former one, we 
prove Theorem~\ref{thm:main:minimal} first.

\begin{proposition} \label{prop:omega:subsoc}
Let $G$ be a minimal abelian group that admits no non-trivial 
quasi-convex null sequences. Then the completion $\widetilde G$ of $G$ 
contains no closed subgroup $H$ such that $\operatorname{soc}(H)$ 
contains  a~non-trivial null sequence that is quasi-convex in $H$.
\end{proposition}

\begin{proof}
Let $\{x_n\} \subseteq\operatorname{soc}(H)$  be a non-trivial null sequence that is quasi-convex in $H$.
Since $H$ is closed in $\widetilde G $, by Corollary~\ref{cor:prel:subgrpB}, $\{x_n\}$ is quasi-convex in
$\widetilde G $. By Corollary~\ref{cor:prel:st}(a),
$\operatorname{soc}(\widetilde G)\subseteq G $, and thus $\{x_n\} \subseteq G $. 
Therefore, by Corollary~\ref{cor:prel:subgrpA}, $\{x_n\}$ is a non-trivial 
quasi-convex null sequence in~$G$, contrary to our assumption.
\end{proof}

\begin{proposition} \label{prop:omega:noTJZ}
Let $G$ be a minimal abelian group that admits no non-trivial quasi-convex 
null sequences. Then the completion $\widetilde G$ of $G$ contains no 
subgroups that are topologically isomorphic to:

\begin{myalphlist}

\item
$\mathbb{J}_p$ \ for some prime $p$,

\item
$\mathbb{T}$, or

\item
$\prod\limits_{k=1}^\infty \mathbb{Z}_{r_k}$ \ for square-free numbers 
$r_k > 3$.
\end{myalphlist}
\end{proposition}

\begin{proof}

(a) Assume that $H$ is a subgroup of $\widetilde G$ that is topologically 
isomorphic to $\mathbb{J}_p$ for some prime $p$. By 
Theorem~\ref{thm:prel:ess}(a), there is $y\in G \cap H$ such that $y \neq 0$. 
Since $\operatorname{cl}_{\widetilde G} \langle y\rangle\cong \mathbb{J}_p$,
by replacing $H$ with $\operatorname{cl}_{\widetilde G} \langle y\rangle$ 
if necessary, we may identify $\langle y \rangle$ with $\mathbb{Z}$ in 
$\mathbb{J}_p$. Thus, by Example~\ref{exp:intro:seqs}(a), $H$ admits a 
non-trivial quasi-convex sequence $\{x_n\}$ such that 
$\{x_n\} \subseteq G \cap H $. Since $H$ is closed in $\widetilde G$, by 
Corollary~\ref{cor:prel:subgrpB}, $\{x_n\}$ is quasi-convex in $\widetilde 
G$. Therefore, by Corollary~\ref{cor:prel:subgrpA}, $\{x_n\}$ is a 
non-trivial quasi-convex null sequence in $G$, contrary to our assumption.

(b) Assume that $H$ is a subgroup of $\widetilde G$ that is topologically 
isomorphic to $\mathbb{T}$. Then $H$ is closed in $\widetilde G$, and by 
Example~\ref{ex:T:socT}, $H$ admits a non-trivial quasi-convex null 
sequence $\{x_n\}$ such that $\{x_n\} \subseteq \operatorname{soc}(H)  
\cong \operatorname{soc}(\mathbb{T})$. By 
Proposition~\ref{prop:omega:subsoc}, the statement follows.

(c) Assume that $H$ is a subgroup of $\widetilde G$ that is topologically 
isomorphic to $\prod\limits_{k=1}^\infty \mathbb{Z}_{r_k}$, where $r_k > 3$ 
and $r_k$ is square-free for every $k$. Then $H$ is 
closed in $\widetilde G $, and by Example~\ref{exp:intro:seqs}(c), $H$ 
admits a~non-trivial quasi-convex null sequence $\{x_n\}$ such that 
$\{x_n\} \subseteq \operatorname{soc}(H) \cong 
\bigoplus\limits_{k=1}^\infty \mathbb{Z}_{r_k}$. The statement follows now 
by Proposition~\ref{prop:omega:subsoc}.
\end{proof}

\begin{lemma} \label{lemma:omega:noTJZ}
Let $K$ be a compact abelian group that contains no subgroups that are topologically isomorphic to $\mathbb{T}$, $\mathbb{J}_p$ for some $p$, or $\prod\limits_{k=1}^\infty \mathbb{Z}_{r_k}$ \ for square-free numbers $r_k  >  3$. Then $K  =  K_p \times  F $, where  $K_p$ is a compact bounded abelian $p$-group, $p  \leq  3$, and $F$ is a finite abelian group.
\end{lemma}

\begin{proof} {\itshape Step 1.} Suppose that $K$ is a pro-finite group. Since $K$ contains no subgroups that are topologically isomorphic to 
$\mathbb{J}_p$ for some prime $p$, it is an exotic torus. The group $K$ has no connected quotients, because it is pro-finite, and thus, by Theorem~\ref{thm:prel:ET}, $K= \prod\limits_{p} K_p$, where each $K_p$ is a~compact bounded abelian $p$-group. Consequently, each
$K_p$ is topologically isomorphic to a~product of finite cyclic $p$-groups (cf.~\cite[4.2.2]{DikProSto}), and $K_p$ is infinite if and only 
if it contains a~subgroup that is topologically isomorphic to $\mathbb{Z}_p^\omega$. By our assumption, $K$  contains no such subgroups for $p>3$. Hence, $K_p$ is finite for $p > 3$.

Put $K^\prime  := \prod\limits_{p>3} K_p$. If $K^\prime$ is infinite, then  there are infinitely many primes $p_k > 3$ such that $K_{p_k} \neq   0$. Consequently, $K^\prime$ (and thus $K$) contains a  subgroup that is topologically isomorphic to the product $\prod\limits_{k=1}^\infty \mathbb{Z}_{p_k}$, contrary to our assumption. This shows that $K^\prime$ is finite. 

Finally, if both $K_2$ and $K_3$ are infinite, then $K$ contains a 
subgroup that is topologically isomorphic to $\mathbb{Z}_2^\omega \times 
\mathbb{Z}_3^\omega \cong \mathbb{Z}_6^\omega$, contrary to our 
assumption. Hence, one of $K_2$ and $K_3$ is finite,~and either $K = K_2 
\times F$, where $F := K_3 \times K^\prime$ is finite, or $K = K_3 \times 
F $, where $F := K_2 \times K^\prime$ is finite.

{\itshape Step 2.}
In the general case, we show that $K$ is pro-finite. To that end, let $C$ 
be the connected component of $K$. Since $K$ is an exotic torus, so is 
$C$, and by Theorem~\ref{thm:prel:ET}, $C$ contains a~closed subgroup $B$ 
such that $B = \prod\limits_{p} B_p$, where each $B_p$ is a finite 
$p$-group, and $C/B \cong \mathbb{T}^n$ for some $n \in \mathbb{N}$. In 
particular, $B$~is a compact pro-finite group that satisfies the 
conditions of this lemma. Thus, by what we have shown so far, $B^\prime:= 
\prod\limits_{p>3} B_p$ is finite, and therefore $B = B_2 \times B_3 
\times B^\prime$ is finite. Consequently, by Pontryagin duality, $\widehat 
B \cong \widehat C / B^\perp$ is finite (cf.~\cite[Theorem~54]{Pontr}), 
and $B^\perp \cong \widehat{C/B} = \mathbb{Z}^n$ 
(cf.~\cite[Theorem~37]{Pontr}). This implies that $\widehat C$ is finitely 
generated. On the other hand, $\widehat C$ is torsion free, because $C$ is 
connected (cf.~\cite[Example~73]{Pontr}), which means that 
$\widehat C = \mathbb{Z}^n$ and $C \cong \mathbb{T}^n$. By our  
assumption, however, $K$ contains no subgroup that is topologically 
isomorphic to $\mathbb{T}$.  Hence, $n=0$ and $C=0$. This shows that $K= 
B$ is pro-finite, and the statement follows from Step 1.
\end{proof}

Proposition~\ref{prop:omega:noTJZ} combined with 
Lemma~\ref{lemma:omega:noTJZ} yields the following consequence.

\begin{corollary} \label{cor:omega:split}
Let $G$ be a minimal abelian group that admits no non-trivial quasi-convex
null sequences. Then the completion $\widetilde G$ of $G$ is a bounded 
compact abelian group, and $\widetilde G =  K_p  \times F$,  where  $K_p$ 
is a compact bounded abelian  $p$-group, $p\leq 3$, and $F$ is a finite 
abelian group.
\qed
\end{corollary}

\begin{proof}[Proof of Theorem~\ref{thm:main:minimal}.]
(i) $\Rightarrow$ (ii):
Let $G$ be a minimal abelian group that admits no non-trivial quasi-convex 
null sequences. By Corollary~\ref{cor:omega:split}, $\widetilde G = K_p 
\times F$, where $K_p$ is a compact bounded abelian $p$-group ($p \leq 
3$), and $F$ is a finite abelian group. Without loss of generality, we may 
assume that $F$ contains no $p$-elements. Let $e=p^a m$ be the exponent 
of~$\widetilde G$,~where $m$ and $p$ are coprime. Then $p^aG$ is dense in 
$p^a\widetilde G = F$, and thus $F = p^a G \subseteq G $. Therefore, for 
$P := G \cap K_p$, one has $G = P \times F$, and $P$ is a bounded abelian 
$p$-group. Since $P$ is a closed subgroup of $G$, by 
Corollary~\ref{cor:prel:subgrpB}, every quasi-convex sequence in $P$ is 
also quasi-convex in $G$, and so $P$ admits no non-trivial quasi-convex 
null sequences. In order to show that $P$ is minimal, let $H$ be a closed 
subgroup of $\widetilde P \subseteq K_p$. Then, in particular, $H$ is a 
closed subgroup of $\widetilde G $. Consequently, by 
Theorem~\ref{thm:prel:ess}(a),
\begin{align}
P\cap H = (G\cap K_p)\cap H = G \cap H \neq \{0\},
\end{align}
because $G$ is minimal. Hence, by Theorem~\ref{thm:prel:ess}(a), $P$ is 
minimal. (It is well known that closed central subgroups of minimal groups 
are minimal, but in this case, a direct proof was also available.)

(ii) $\Rightarrow$ (iii):  By Theorem~\ref{thm:prel:ess}(a), $P$ is 
precompact. Thus, by Theorem~\ref{thm:main:BP}, $P[2]$ or $P[3]$ is 
sequentially open in $P$. Since $F$ is finite, $P$ is open in $G$. 
Therefore, $P[2]$ or $P[3]$ is sequentially open in~$G$. Hence,
$G[2]$ or $G[3]$ is sequentially open in $G$.

(iii) $\Rightarrow$ (iv): Since $G$ is minimal, by 
Corollary~\ref{cor:prel:st}(a), $G[p]=\widetilde G[p]$ for every prime $p$, 
and in particular $G[p]$ is compact. Thus, 
$G[p]\cong \mathbb{Z}_p^{\kappa_p}$ for some cardinal $\kappa_p$ for every 
prime $p$.

(iv) $\Rightarrow$ (i): If $G$ contains a sequentially open compact 
subgroup of the form $\mathbb{Z}_2^\kappa$ or $\mathbb{Z}_3^\kappa$, then 
$G[2]$  or $G[3]$ is sequentially open in $G$, and therefore the statement 
follows by Proposition~\ref{prop:bad:G23}.
\end{proof}

We turn now to the  case where $G$ is $\omega$-bounded.

\begin{proposition} \label{prop:omega:oB}
Let $G$ be an $\omega$-bounded abelian group that admits no non-trivial quasi-convex null sequences. Then:

\begin{myalphlist}

\item
$G$ is bounded;

\item
the subgroup $G_p$  of $p$-elements is finite for 
$p>3$;

\item
at least one of $G_2$ and $G_3$ is finite;

\item
$2G_2$ and $3G_3$ are finite.

\end{myalphlist}
\end{proposition}

The proof of Proposition~\ref{prop:omega:oB} is based on a
well-known result of Comfort and Robertson.

\begin{ftheorem}[{\cite[7.4]{ComfRob}}] \label{thm:omega:psATB}
Every pseudocompact abelian torsion group is bounded.
\end{ftheorem}

\begin{proof}[Proof of Proposition~\ref{prop:omega:oB}.]
(a) Let $x\in G$. Since $G$ is $\omega$-bounded, $\langle x \rangle$ is 
contained in a compact subgroup of $G$, and thus 
$K := \operatorname{cl}_G \langle x \rangle$ is compact. If $x$ has an 
infinite order, then $2K$ and $3K$ are infinite, and so by  
Theorem~\ref{thm:intro:LCAqcs}, $K$ admits a non-trivial  quasi-convex 
null sequence $\{x_n\}$. Since $K$ is a closed subgroup of $G$, by 
Corollary~\ref{cor:prel:subgrpB}, $\{x_n\}$ is a non-trivial quasi-convex 
null sequence in $G$, contrary to our assumption. This shows that every  
element in $G$ has a finite order; in other words, $G$ is a torsion group.
Since $G$ is $\omega$-bounded, in particular, it is pseudocompact. Hence, 
by Theorem~\ref{thm:omega:psATB}, $G$ is bounded.

(b) Let $p$ be a prime, and suppose that $G_p$ is infinite. Then $G[p]$ 
contains a countably infinite subset $S$, which in turn is contained in a 
compact subgroup $K$, because $G$ is  $\omega$-bounded. Since $K[p]$ is 
an infinite compact group of exponent $p$, it is topologically isomorphic  
to $\mathbb{Z}_p^\lambda$ for some infinite cardinal $\lambda$ 
(cf.~\cite[4.2.2]{DikProSto}). In particular, $G$ contains a closed 
subgroup $H$ that is topologically isomorphic to $\mathbb{Z}_p^\omega$.

If $p> 3$, then by Example~\ref{exp:intro:seqs}(c), $H$ admits a 
non-trivial quasi-convex null sequence, and since $H$ is closed in $G$, 
this sequence will also be quasi-convex in $G$ according to 
Corollary~\ref{cor:prel:subgrpB}, contrary to our assumption. This shows 
that $G_p$ is finite for $p>3$.

(c) Assume that both $G_2$ and $G_3$ are infinite. Then,  by what we have seen so far, $G$ contains a subgroup $H_2$ that is topologically 
isomorphic to $\mathbb{Z}_2^\omega$, and a subgroup $H_3$ that is topologically isomorphic to $\mathbb{Z}_3^\omega$. Thus,
$H  := H_2 + H_3$ is topologically isomorphic to $\mathbb{Z}_6^\omega$. By Example~\ref{exp:intro:seqs}(c), $H$~admits a non-trivial quasi-convex null sequence, and since $H$ is closed in $G$, this sequence will also be quasi-convex in $G$ according to Corollary~\ref{cor:prel:subgrpB}, contrary to our assumption. This shows that at least one of $G_2$ and $G_3$ is finite.

(d) Let $p = 2$ or $p = 3$, and assume that $pG_p$ is infinite. Then, in particular, $(pG_p)[p]$ is infinite, and so there is a countably infinite subset  $S$ of $G$ such that $pS$ is infinite and $p^2 S  =  0$. Since~$G$ is  $\omega$-bounded, $S$ is contained in a compact subgroup $K$ of 
$G$. By replacing $K$ with $K[p^2]$ if necessary, we may assume that $K$ has exponent $p^2$, and so $K$ is topologically isomorphic to
$\mathbb{Z}_p^{\lambda_1} \times \mathbb{Z}_{p^2}^{\lambda_2}$ for~some cardinals $\lambda_1$ and $\lambda_2$. As 
$pS \subseteq  pK  \cong \mathbb{Z}_{p}^{\lambda_2}$ is infinite, $\lambda_2$ is infinite, and thus $G$ contains a subgroup $H$ 
that is topologically isomorphic to $\mathbb{Z}_{p^2}^\omega$. By Example~\ref{exp:intro:seqs}(c), $H$ admits a non-trivial 
quasi-convex null sequence, and since $H$ is closed in $G$, this sequence will also be quasi-convex in $G$ according to
Corollary~\ref{cor:prel:subgrpB}, contrary to our assumption. This shows that $2G_2$ and $3G_3$ are finite.
\end{proof}

We are now ready to prove  Theorem~\ref{thm:main:omega}.

\begin{proof}[Proof of Theorem~\ref{thm:main:omega}.]
By Proposition~\ref{prop:bad:G23}, the implication (ii) $\Rightarrow$ (i)  
holds for every precompact group, and obviously so does the equivalence 
(ii) $\Leftrightarrow$ (iii). The implication (iv) $\Rightarrow$ (iii) is 
also clear. Thus, it suffices to prove that (i) $\Rightarrow$ (iii), and 
if $G$ is totally minimal, (i) $\Rightarrow$ (iv).

Suppose that $G$ is $\omega$-bounded and admits no non-trivial 
quasi-convex null sequences.  Then, by Proposition~\ref{prop:omega:oB}(a), 
$G$ is bounded, and so $\widetilde G$ is bounded. Thus, $\widetilde G$ is 
a product of its subgroups of $p\mbox{-}$elements, and therefore 
$G\cong  G_2 \times G_3 \times G_{p_1}  \times   \cdots   \times G_{p_l}$, 
where  $p_1,\ldots,p_l  >  3$ are 
prime factors of the exponent of $G $. By 
Proposition~\ref{prop:omega:oB}(b), each $G_{p_i}$ is finite. By 
Proposition~\ref{prop:omega:oB}(c), one of $G_2$ and $G_3$ is finite, and so 
$G  \cong   G_p  \times   F$, where 
$p = 2$ or $p = 3$, and $F$ is a finite abelian group of order coprime to 
$p$. By Proposition~\ref{prop:omega:oB}(d), $pG_p$ is finite, and hence 
$pG \cong pG_p \times F$ is finite, as desired.

Suppose that $G$ is totally minimal and admits no non-trivial quasi-convex 
null sequences. Then, by Corollary~\ref{cor:omega:split}, the completion 
$\widetilde G$ of $G$ is bounded. Thus, by Corollary~\ref{cor:prel:st}(b), 
$\widetilde G = G$, and so $G$ is compact. Hence, both (iii) and (iv) 
follow from Theorem~\ref{thm:intro:LCAqcs}.
\end{proof}

\section*{Acknowledgements}

The authors are grateful to Karen Kipper for her kind help in proofreading 
this paper for grammar and punctuation. The authors wish to thank the 
anonymous referee for the wealth of constructive comments that led to an 
improved presentation of this paper.

{\footnotesize

\bibliography{notes,notes2,notes3}
}

\begin{samepage}

\bigskip
\noindent
\begin{tabular}{l @{\hspace{1.8cm}} l}
Department of Mathematics and Computer Science & Halifax, Nova Scotia\\
University of Udine & Canada\\
Via delle Scienze, 208 -- Loc. Rizzi, 33100 Udine &  \\
Italy &  \\ 
& \\
\em email: dikranja@dimi.uniud.it  &
\em email: lukacs@topgroups.ca
\end{tabular}

\end{samepage}

\end{document}